\documentclass[12pt]{article}

\usepackage{amsmath,amssymb}

\begin{document}

\title{Euler-Poincar\'e functions}
\author{Anton Deitmar}

\date{}
\maketitle

\tableofcontents

\setlength{\parskip}{7pt}

\def \1{{\bf 1}}
\def \a{{{\mathfrak a}}}
\def \ad{{\rm ad}}
\def \al{\alpha}
\def \ar{{\alpha_r}}
\def \A{{\mathbb A}}
\def \Ad{{\rm Ad}}
\def \Aut{{\rm Aut}}
\def \b{{{\mathfrak b}}}
\def \bs{\backslash}
\def \B{{\cal B}}
\def \c{{\mathfrak c}}
\def \cent{{\rm cent}}
\def \C{{\mathbb C}}
\def \CA{{\cal A}}
\def \CB{{\cal B}}
\def \CC{{\cal C}}
\def \CE{{\cal E}}
\def \CF{{\cal F}}
\def \CG{{\cal G}}
\def \CH{{\cal H}}
\def \CHC{{\cal HC}}
\def \CL{{\cal L}}
\def \CM{{\cal M}}
\def \CN{{\cal N}}
\def \CP{{\cal P}}
\def \CQ{{\cal Q}}
\def \CO{{\cal O}}
\def \CS{{\cal S}}
\def \CT{{\cal T}}
\def \CV{{\cal V}}
\def \det{{\rm det}}
\def \diag{{\rm diag}}
\def \dist{{\rm dist}}
\def \End{{\rm End}}
\def \Fx{{\mathfrak x}}
\def \FX{{\mathfrak X}}
\def \g{{{\mathfrak g}}}
\def \ga{\gamma}
\def \Ga{\Gamma}
\def \h{{{\mathfrak h}}}
\def \Hom{{\rm Hom}}
\def \im{{\rm im}}
\def \Im{{\rm Im}}
\def \Ind{{\rm Ind}}
\def \k{{{\mathfrak k}}}
\def \K{{\cal K}}
\def \l{{\mathfrak l}}
\def \la{\lambda}
\def \lap{\triangle}
\def \La{\Lambda}
\def \m{{{\mathfrak m}}}
\def \mod{{\rm mod}}
\def \n{{{\mathfrak n}}}
\def \name{\bf}
\def \N{\mathbb N}
\def \o{{\mathfrak o}}
\def \ord{{\rm ord}}
\def \O{{\cal O}}
\def \p{{{\mathfrak p}}}
\def \ph{\varphi}
\def \prf{\noindent{\bf Proof: }}
\def \Per{{\rm Per}}
\def \q{{\mathfrak q}}
\def \qed{\ifmmode\eqno Q.E.D.\else\noproof\vskip 12pt plus 3pt minus 9pt \fi}
 \def\noproof{{\unskip\nobreak\hfill\penalty50\hskip2em\hbox{}%
     \nobreak\hfill Q.E.D.\parfillskip=0pt%
     \finalhyphendemerits=0\par}}
\def \Q{\mathbb Q}
\def \res{{\rm res}}
\def \R{{\mathbb R}}
\def \Re{{\rm Re \hspace{1pt}}}
\def \r{{\mathfrak r}}
\def \ra{\rightarrow}
\def \rank{{\rm rank}}
\def \Rep{{\rm Rep}}
\def \supp{{\rm supp}}
\def \Spin{{\rm Spin}}
\def \t{{{\mathfrak t}}}
\def \T{{\mathbb T}}
\def \tr{{\hspace{1pt}\rm tr\hspace{2pt}}}
\def \vol{{\rm vol}}
\def \z{\zeta}
\def \Z{\mathbb Z}
\def \={\ =\ }
\newcommand{\choice}[4]{\left\{
            \begin{array}{cl}#1&#2\\ #3&#4\end{array}\right.}
\newcommand{\rez}[1]{\frac{1}{#1}}
\newcommand{\der}[1]{\frac{\partial}{\partial #1}}
\newcommand{\norm}[1]{\parallel #1 \parallel}
\renewcommand{\matrix}[4]{\left( \begin{array}{cc}#1 & #2 \\ #3 & #4 \end{array}
            \right)}

\newcounter{lemma}
\newcounter{corollary}
\newcounter{proposition}
\newcounter{theorex}

\newtheorem{theorem}{Theorem}[section]
\newtheorem{conjecture}[theorem]{Conjecture}
\newtheorem{lemma}[theorem]{Lemma}
\newtheorem{corollary}[theorem]{Corollary}
\newtheorem{proposition}[theorem]{Proposition}

\begin{center} {\bf Introduction} \end{center}

Euler-Poincar\'e functions and pseudo-coefficients are
important tools in harmonic analysis since they help to
single out a given representation in direct integrals or
sums. Given a real reductive group $G$ that admits a
compact Cartan subgroup, these functions can be attached
to a given finite dimensional representation $\tau$ of a
maximal compact subgroup $K$. In \cite{Lab} a construction
was given in the case when $\tau$ extends to the bigger
group $G$. In this paper we give a new construction for an
arbitrary representation $\tau$. Instead, we put a
condition on the group $G$, namely, that it acts
orientation preservingly on the symmetric space $G/K$.
This condition is satisfied if $G$ is connected or if it
respects a complex structure on $G/K$. For connected $G$
the existence of Euler-Poincar\'e functions is known
\cite{CloDel}, but this is not always sufficient for
applications since Levi components are not in general
connected, even if the ambient group is. After having
constructed Euler-Poincar\'e functions we compute the
orbital integrals of semisimple elements.

$ $

\section{Notations} \label{notations}

We denote Lie groups by upper case roman letters
$G,H,K,\dots$ and the corresponding real Lie algebras by
lower case German letters with index $0$, that is:
$\g_0,\h_0,\k_0,\dots$. The complexified Lie algebras will
be denoted by $\g,\h,\k,\dots$, so, for example: $\g
=\g_0\otimes_\R \C$.

Let $G$ be a real reductive group \cite{wall-rg1} and fix
a maximal compact subgroup $K$ of $G$. Fix a Cartan
involution $\theta$ with fixed point set the maximal
compact subgroup $K$ and let $\k_0$ be the Lie algebra of
$K$.  The group $K$ acts on $\g_0$ via the adjoint
representation and there is a $K$-stable decomposition
$\g_0=\k_0\oplus\p_0$,  where $\p_0$ is the eigenspace of
(the differential of) $\theta$ to the eigenvalue $-1$.
Write $\g=\k\oplus\p$ for the complexification. This is
called the {\it Cartan decomposition}.

Fix a symmetric bilinear form $B : \g_0\times\g_0\ra \R$
such that
\begin{itemize}
\item
$B$ is invariant, that is $B(\Ad(g)X,\Ad(g)Y)=B(X,Y)$ for
all $g\in G$ and all $X,Y\in\g_0$ and
\item
$B$ is negative definite on $\k_0$ and positive definite
on its orthocomplement $\p_0=\k_0^\perp\subset\g_0$.
\end{itemize}
When $G$ is semisimple we can choose $B$ to be the Killing
form of $\g_0$.

Let $U(\g)$ be the universal enveloping algebra of $\g$,
then $B$ gives rise to a Casimir element $C\in U(\g)$. For
an irreducible admissible representation $(\pi, V_\pi)$
the operator $C$ will act on $V_\pi$ by a scalar
$\pi(C)\in\C$.

Let $X$ denote the quotient manifold $G/K$. The tangent
space at $eK$ identifies with $\p_0$ and the form $B$
gives a $K$-invariant positive definite inner product on
this space. Translating this by elements of $G$ defines a
$G$-invariant Riemannian metric on $X$. This makes $X$ the
most general globally symmetric space of the noncompact
type \cite{helg}.

Let $\hat{G}$\index{$\hat{G}$} denote the {\it unitary
dual}\index{unitary dual} of $G$, i.e., $\hat{G}$ is the
set of isomorphism classes of irreducible unitary
representations of $G$.

The form $\langle X,Y\rangle =-B(X,\theta(Y)$ is positive
definite on $\g_0$ and therefore induces a positive
definite left invariant top differential form $\omega_L$
on any closed subgroup $L$ of $G$. If $L$ is compact we set
$$
v(L)\= \int_L\omega_L.
$$
Let $H=AB$ be a $\theta$-stable Cartan subgroup where $A$
is the connected split component of $H$ and $B$ is
compact. The double use of the letter $B$ here will not
cause any confusion. Then $B\subset K$. Let $\Phi$ denote
the root system of $(\g,\h)$, where $\g$ and $\h$ are the
complexified Lie algebras of $G$ and $H$. Let $\g
=\h\oplus\bigoplus_{\alpha\in\Phi}\g_\alpha$ be the root
space decomposition. Let $x\ra x^c$ denote the complex
conjugation on $\g$ with respect to the real form $\g_0=
Lie(G)$. A root $\alpha$ is called {\it
imaginary}\index{imaginary root} if $\alpha^c=-\alpha$ and
it is called {\it real}\index{real root} if
$\alpha^c=\alpha$. Every root space $\g_\alpha$ is one
dimensional and has a generator $X_\alpha$ satisfying:
$$
[X_\alpha,X_{-\alpha}]=Y_\alpha\ \ \ \ {\rm with}\ \ \ \
\alpha(.)=B(Y_\alpha,.)
$$ $$
B(X_\alpha,X_{-\alpha})=1
$$
and $X_\alpha^c=X_{\alpha^c}$ if $\alpha$ is non-imaginary
and $X_\alpha^c=\pm X_{-\alpha}$ if $\alpha$ is imaginary.
An imaginary root $\alpha$ is called {\it
compact}\index{compact root} if $X_\alpha^c=-X_{-\alpha}$
and {\it noncompact} otherwise. Let $\Phi_n$ be the set of
noncompact imaginary roots and choose a set $\Phi^+$ of
positive roots such that for $\alpha\in\Phi^+$
nonimaginary we have that $\alpha^c\in\Phi^+$. Let
$W=W(G,H)$ be the {\it Weyl group} of $(G,H)$\index{Weyl
group of $(G,H)$}, that is
$$
W\= \frac{{\rm normalizer}(H)}{{\rm centralizer}(H)}.
$$
Let ${\rm rk}_\R(G)$ be the dimension of a maximal
$\R$-split torus in $G$ and let $\nu =\dim G/K -{\rm
rk}_\R(G)$. We define the {\it Harish-Chandra constant} of
$G$ by
$$
c_G\= (-1)^{|\Phi_n^+|} (2\pi)^{|\Phi^+|}
2^{\nu/2}\frac{v(T)}{v(K)}|W|.
$$

\section{Normalization of Haar measures}
Although the results will not depend on normalizations we
will need to normalize Haar measures for the computations
along the way. First for any compact subgroup $C\subset G$
we normalize its Haar measure so that it has total mass
one, i.e., $\vol(C)=1$. Next let $H\subset G$ be a
reductive subgroup, and let $\theta_H$ be a Cartan
involution on $H$ with fixed point set $K_H$. The same way
as for $G$ itself the form $B$ restricted to the Lie
algebra of $H$ induces a Riemannian metric on the manifold
$X_H=H/K_H$. Let $dx$ denote the volume element of that
metric. We get a Haar measure on $H$ by defining
$$
\int_Hf(h) dh \= \int_{X_H}\int_{K_H} f(xk) dk dx
$$
for any continuous function of compact support $f$ on $H$.

Let $P\subset G$ be a {\it parabolic
subgroup}\index{parabolic subgroup} of $G$ (see
\cite{wall-rg1} 2.2). Let $P=MAN$ be the {\it Langlands
decomposition} \index{Langlands decomposition} of $P$.
Then $M$ and $A$ are reductive, so there Haar measures can
be normalized as above. Since $G=PK=MANK$ there is a
unique Haar measure $dn$ on the unipotent radical $N$ such
that for any constant function $f$ of compact support on
$G$ it holds:
$$
\int_G f(x)dx \= \int_M\int_A\int_N\int_K f(mank) dkdndadm.
$$
Note that these normalizations coincide for Levi subgroups
with the ones met by Harish-Chandra in (\cite{HC-HA1}
sect. 7).

\section{Invariant distributions}
In this section we shall throughout assume that $G$ is a
real reductive group of inner type \cite{wall-rg1}. A
distribution $T$ on $G$, i.e., a continuous linear
functional $T: C_c^\infty(G)\ra \C$ is called {\it
invariant}\index{invariant distribution} if for any $f\in
C_c^\infty(G)$ and any $y\in G$ it holds: $T(f^y)=T(f)$,
where $f^y(x)=f(yxy^{-1})$. Examples are:
\begin{itemize}
\item orbital integrals: $f\mapsto \CO_g(f)=\int_{G_g\bs G}f(x^{-1}gx)dx$ and
\item traces: $f\mapsto \tr\pi(f)$ for $\pi\in\hat{G}$.
\end{itemize}
These two examples can each be expressed in terms of the
other. Firstly, Harish-Chandra proved that for any
$\pi\in\hat{G}$ there exists a conjugation invariant
locally integrable function $\Theta_\pi$ on $G$ such that
for any $f\in C_c^\infty(G)$
$$
\tr\pi(f)\=\int_Gf(x)\Theta_\pi(x)  dx.
$$
Recall the {\it Weyl integration formula} which says that
for any integrable function $\ph$ on $G$ we have
$$
\int_G\ph(x)dx\=
\sum_{j=1}^r\rez{|W(G,H_j)|}\int_{H_j}\int_{G/H_j}\ph(xhx^{-1})|\det(1-h|\g
/\h_j)|dx dh,
$$
where $H_1,\dots,H_r$ is a maximal set of nonconjugate
Cartan subgroups in $G$ and for each Cartan subgroup $H$
we let $W(G,H)$\index{$W(G,H)$} denote its Weyl group,
i.e., the quotient of the normalizer of $H$ in $G$ by its
centralizer.

An element $x$ of $G$ is called {\it
regular}\index{regular} if its centralizer is a Cartan
subgroup. The set of regular elements $G^{reg}$ is open
and dense in $G$ therefore the integral above can be taken
over $G^{reg}$ only. Letting $H_j^{reg}:=H_j\cap G^{reg}$
we get

\begin{proposition}
Let $N$ be a natural number bigger than $\frac{\dim
G}{2}$, then for any $f\in L_{2N}^1(G)$ and any
$\pi\in\hat{G}$ we have
$$
\tr\pi(f)\=\sum_{j=1}^r\rez{|W(G,H_j)|}\int_{H_j^{reg}}
\CO_h(f)\Theta_\pi(h)|\det(1-h|\g /\h_j)| dh.
$$
\end{proposition}

That is, we have expressed the trace distribution in terms
of orbital integrals. In the other direction it is also
possible to express semisimple orbital integrals in terms
of traces.

At first let $H$ be a $\theta$-stable Cartan subgroup of
$G$. Let $\h$ be its complex Lie algebra and let $\Phi
=\Phi(\g,\h)$ be the set of roots. Let $x\ra x^c$ denote
the complex conjugation on $\g$ with respect to the real
form $\g_0= Lie_\R(G)$. Choose an ordering $\Phi^+\subset
\Phi$ and let $\Phi^+_I$ be the set of positive imaginary
roots. To any root $\alpha\in\Phi$ let
\begin{eqnarray*}
H &\ra & \C^\times\\
 h &\mapsto & h^\alpha
\end{eqnarray*}\index{$h^\alpha$}
be its character, that is, for $X\in\g_\alpha$ the root
space to $\alpha$ and any $h\in H$ we have
$\Ad(h)X=h^\alpha X$. Now put
$$
\ '\lap_I(h)\= \prod_{\alpha\in\Phi_I^+}(1-h^{-\alpha}).
$$\index{$\ '\lap_I(h)$}
Let $H=AT$ where $A$ is the connected split component and
$T$ is compact. An element $at\in AT=H$ is called {\it
split regular}\index{split regular} if the centralizer of
$a$ in $G$ equals the centralizer of $A$ in $G$. The split
regular elements form a dense open subset containing the
regular elements of $H$. Choose a parabolic $P$ with split
component $A$, so $P$ has Langlands decomposition $P=MAN$.
For $at\in AT =H$ let
\begin{eqnarray*}
\lap_+(at) &=& \left| \det((1-\Ad((at)^{-1}))|_{\g
/\a\oplus
\m})\right|^{\rez{2}}\\ {}\\
 & =& \left|\det((1-\Ad((at)^{-1}))|_\n)\right| a^{\rho_P}\\ {}\\
 &=&
 \left|\prod_{\alpha\in\Phi^+-\Phi_I^+}(1-(at)^{-\alpha})\right|
 a^{\rho_P},
\end{eqnarray*}\index{$\lap_+(at)$}
where $\rho_P$\index{$\rho_P$} is the half of the sum of
the roots in $\Phi(P,A)$, i.e., $a^{2\rho_P}=\det(a|\n)$.
We will also write $h^{\rho_P}$ instead of $a^{\rho_P}$.

For any $h\in H^{reg} = H\cap G^{reg}$ let
$$
'F_f^H(h)\= 'F_f(h)\= \ '\lap_I(h) \lap_+(h) \int_{G/A}
f(xhx^{-1}) dx.
$$
It then follows directly from the definitions that for
$h\in H^{reg}$ it holds
$$
\CO_h(f) \= \frac{'F_f(h)}
                 {h^{\rho_P}\det(1-h^{-1}|(\g/\h)^+)},
$$
where $(\g/\h)^+$ is the sum of the root spaces attached to
positive roots. There is an extension of this identity to
nonregular elements as follows: For $h\in H$ let $G_h$
denote its centralizer in $G$. Let $\Phi^+(\g_h,\h)$ be
the set of positive roots of $(\g_h,\h)$. Let
$$
\varpi_h \= \prod_{\alpha\in\Phi^+(\g_h,\h)}Y_\alpha,
$$
then $\varpi_h$ defines a left invariant differential
operator on $G$.

\begin{lemma}
For any $f\in L_{2N}^1(G)$ and any $h\in H$ we have
$$
\CO_h(f) \= \frac{\varpi_h 'F_f(h)}
                 {c_h h^{\rho_P}\det(1-h^{-1}|(\g/\g_h)^+)}.
$$
\end{lemma}

\prf This is proven in section 17 of \cite{HC-DS}. \qed

Our aim is to express orbital integrals in terms of traces
of representations. By the above lemma it is enough to
express $'F_f(h)$ it terms of traces of $f$ when $h\in
H^{reg}$. For this let $H_1=A_1T_1$ be another
$\theta$-stable Cartan subgroup of $G$ and let
$P_1=M_1A_1N_1$ be a parabolic with split component $A_1$.
Let $K_1=K\cap M_1$. Since $G$ is connected the compact
group $T_1$ is an abelian torus and its unitary dual
$\widehat{T_1}$ is a lattice. The Weyl group
$W=W(M_1,T_1)$ acts on $\widehat{T_1}$ and
$\widehat{t_1}\in\widehat{T_1}$ is called {\it regular} if
its stabilizer $W(\widehat{t_1})$ in $W$ is trivial. The
regular set $\widehat{T_1}^{reg}$ modulo the action of
$W(K_1,T_1)\subset W(M_1,T_1)$ parameterizes the discrete
series representations of $M_1$ (see \cite{Knapp}). For
$\widehat{t_1}\in\widehat{T_1}$ Harish-Chandra \cite{HC-S}
defined a distribution $\Theta_{\widehat{t_1}}$ on $G$
which happens to be the trace of the discrete series
representation $\pi_{\widehat{t_1}}$ attached to
$\widehat{t_1}$ when $\widehat{t_1}$ is regular. When
$\widehat{t_1}$ is not regular the distribution
$\Theta_{\widehat{t_1}}$ can be expressed as a linear
combination of traces as follows. Choose an ordering of
the roots of $(M_1,T_1)$ and let $\Omega$ be the product
of all positive roots. For any $w\in W$ we have $w\Omega =
\epsilon(w)\Omega$ for a homomorphism $\epsilon : W\ra \{
\pm 1\}$. For nonregular $\widehat{t_1}\in\widehat{T_1}$
we get
$\Theta_{\widehat{t_1}}=\rez{|W(\widehat{t_1})|}\sum_{w\in
W(\widehat{t_1})}\epsilon(w)\Theta'_{w,\widehat{t_1}}$,
where $\Theta'_{w,\widehat{t_1}}$ is the character of an
irreducible representation $\pi_{w,\widehat{t_1}}$ called
a limit of discrete series representation. We will write
$\pi_{\widehat{t_1}}$ for the virtual representation
$\rez{|W(\widehat{t_1})|}\sum_{w\in
W_{\widehat{t_1}}}\epsilon(w)\pi_{w,\widehat{t_1}}$.

Let $\nu :a\mapsto a^\nu$ be a unitary character of $A_1$
then $\widehat{h_1}=(\nu,\widehat{t_1})$ is a character of
$H_1=A_1T_1$. Let $\Theta_{\widehat{h_1}}$ be the
character of the representation $\pi_{\widehat{h_1}}$
induced parabolically from $(\nu,\pi_{\widehat{t_1}})$.
Harish-Chandra has proven

\begin{theorem}\label{inv-orb-int}
Let $H_1,\dots,H_r$ be maximal a set of nonconjugate
$\theta$-stable Cartan subgroups. Let $H=H_j$ for some
$j$. Then for each $j$ there exists a continuous function
$\Phi_{H|H_j}$ on $H^{reg}\times \hat{H_j}$ such that for
$h\in H^{reg}$ it holds
$$
'F_f^H(h)\= \sum_{j=1}^r
\int_{\hat{H_j}}\Phi_{H|H_j}(h,\widehat{h_j})\
\tr\pi_{\widehat{h_j}}(f)\ d\widehat{h_j}.
$$
Further $\Phi_{H|H_j}=0$ unless there is $g\in G$ such that
$gAg^{-1}\subset A_1$. Finally for $H_j=H$ the function
can be given explicitly as
\begin{eqnarray*}
\Phi_{H|H}(h,\hat{h}) &=& \rez{|W(G,H)|}\sum_{w\in
W(G,H)}\epsilon(w|T)\langle w\hat{h},h\rangle\\
 &=& \rez{|W(G,H)|} \ '\lap(h)\Theta_{\hat{h}}(h),
\end{eqnarray*}
where $\ '\lap = \lap_+ \lap_I$.
\end{theorem}

\prf \cite{HC-S}. \qed

\section{Existence of Euler-Poincar\'e functions}

Let $G$ be a real reductive group that admits a compact
Cartan. Fix a maximal compact subgroup $K$ of $G$ and a
Cartan $T$ of $G$ which lies inside  $K$. The group $G$ is
called {\it orientation preserving} \index{orientation
preserving} if $G$ acts by orientation preserving
diffeomorphisms on the manifold $X=G/K$. For example, the
group $G=SL_2(\R)$ is orientation preserving but the group
$PGL_2(\R)$ is not. Recall the Cartan decomposition
$\g_0=\k_0\oplus\p_0$. Note that $G$ is orientation
preserving if and only if its maximal compact subgroup $K$
preserves orientations on $\p_0$.

\begin{lemma} \label{orient}
The following holds:
\begin{itemize}
\item
Any connected group is orientation preserving.
\item
If $X$ carries the structure of a complex manifold which is left
stable by $G$, then $G$ is orientation preserving.
\end{itemize}
\end{lemma}

\prf The first is clear. The second follows from the fact that
biholomorphic maps are orientation preserving.
 \qed

Let $\t$\index{$\t$} be the complexified Lie algebra of
the Cartan subgroup $T$. We choose an ordering of the
roots $\Phi(\g ,\t)$\index{$\Phi(\g ,\t)$} of the pair
$(\g ,\t)$ \cite{wall-rg1}. This choice induces a
decomposition $\p = \p_- \oplus \p_+$,\index{$\p_\pm$}
where $\p_\pm$ is the sum of the positive/negative root
spaces which lie in $\p$. As usual denote by
$\rho$\index{$\rho$} the half sum of the positive roots.
The chosen ordering induces an ordering of the {\it compact
roots}\index{compact root} $\Phi(\k ,\t)$ which form a
subset of the set of all roots $\Phi(\g,\t)$. Let $\rho_K$
denote the half sum of the positive compact roots. Recall
that a function $f$ on $G$ is called {\it $K$-central}
\index{$K$-central} if $f(kxk^{-1})=f(x)$ for all $x\in
G$, $k\in K$. For any $K$-representation $(\rho,V)$ let
$V^K$ denote the space of $K$-fixed vectors and let
$(\breve\tau,V_{\breve\tau})$ denote the dual
representation. The restriction from $G$ to $K$ gives a
ring homomorphism of the representation rings:
$$
res_K^G:\Rep(G)\ra\Rep(K).
$$

\begin{theorem}\label{exist-ep}
(Euler-Poincar\'e functions) Let $(\tau ,V_\tau)$ a finite
dimensional representation of $K$. If either $G$ is orientation
preserving or $\tau$ lies in the image of $res^G_K$, then there
is a compactly supported smooth $K$-central function $f_\tau$ on
$G$ such that for every irreducible unitary representation $(\pi
,V_\pi)$ of $G$ we have
$$
\tr\ \pi (f_\tau) \= \sum_{p=0}^{\dim (\p)} (-1)^p \dim (V_\pi
\otimes \wedge^p\p \otimes V_{\breve{\tau}})^K.
$$
We call $f_\tau$ an {\it Euler-Poincar\'e function} for $\tau$.

If, moreover, $K$ leaves invariant the decomposition
$\p=\p_+\oplus\p_-$ then there is a compactly supported smooth
$K$-central function $g_\tau$ on $G$ such that for every
irreducible unitary representation $(\pi ,V_\pi)$ we have
$$
\tr\ \pi (g_\tau) \= \sum_{p=0}^{\dim (\p_-)} (-1)^p \dim (V_\pi
\otimes \wedge^p\p_- \otimes V_{\breve{\tau}})^K.
$$
\end{theorem}

{\it Remark:} If the representation $\tau$ lies in the
image of $res_K^G$ or the group $G$ is connected then the
theorem is well known, \cite{CloDel}, \cite{Lab}. We only
included it for the sake of completeness. We will therefore
concentrate on the proof in the case when $G$ is
orientation preserving.

\vspace{10pt}

\prf  We will concentrate on the case when $G$ is
orientation preserving,.suffice to say that the other case
can be treated similarly. Without loss of generality
assume that $\tau$ is irreducible. Suppose given a
function $f$ which satisfies the claims of the theorem
except that it is not necessarily $K$-central, then the
function
$$
x\mapsto \int_K f(kxk^{-1})dk
$$
will satisfy all claims of the theorem. Thus one only needs to
construct a function having the claimed traces.

If $G$ is orientation preserving the adjoint action gives a
homomorphism $K\ra SO(\p)$. If this homomorphism happens to lift
to the double cover $Spin(\p)$ \cite{LawMich} we let $\tilde{G}=G$
and $\tilde{K}=K$. In the other case we apply the

\begin{lemma}\label{double-cover} If the homomorphism $K\ra
SO(\p)$ does not factor over $Spin(\p)$ then there is a double
covering $\tilde{G}\ra G$ such that with $\tilde{K}$ denoting the
inverse image of $K$ the induced homomorphism $\tilde{K}\ra
SO(\p)$ factors over $Spin(\p)\ra SO(\p)$. Moreover the kernel of
the map $\tilde{G}\ra G$ lies in the center of $\tilde{G}$
\end{lemma}

\prf At first $\tilde{K}$ is given by the pullback diagram:
$$
\begin{array}{ccc}
\tilde{K} & \ra & Spin(\p)\\
\downarrow & {} & \downarrow \alpha\\
K & \begin{array}{c}\Ad\\ \ra\\ {}\end{array} & SO(\p)
\end{array}
$$
that is, $\tilde{K}$ is given as the set of all $(k,g)\in K\times
Spin(\p)$ such that $\Ad(k)=\alpha(g)$. Then $\tilde{K}$ is a
double cover of $K$.

Next we use the fact that $K$ is a retract of $G$ to show that the
covering $\tilde{K}\ra K$ lifts to $G$. Explicitly let
$P=\exp(\p_0)$ then the map $K\times P\ra G, (k,p)\mapsto kp$ is a
diffeomorphism \cite{wall-rg1}. Let $g\mapsto
(\underline{k}(g),\underline{p}(g))$ be its inverse map. We let
$\tilde{G}\= \tilde{K}\times P$ then the covering $\tilde{K}\ra K$
defines a double covering $\beta :\tilde{G}\ra G$. We have to
install a group structure on $\tilde{G}$ which makes $\beta$ a
homomorphism and reduces to the known one on $\tilde{K}$. Now let
$k,k'\in K$ and $p,p'\in P$ then by
$$
k'p'kp\= k'k\ k^{-1}p'kp
$$
it follows that there are unique maps $a_K : P\times P\ra K$ and
$a_P :P\times P\ra P$ such that
\begin{eqnarray*}
\underline{k}(k'p'kp)&=& k'k a_K(k^{-1}p'k,p)\\
\underline{p}(k'p'kp)&=& a_P(k^{-1}p'k,p).
\end{eqnarray*}
Since $P$ is simply connected the map $a_K$ lifts to a map
$\tilde{a}_K: P\times P\ra \tilde{K}$. Since $P$ is connected
there is exactly one such lifting with $\tilde{a}_K(1,1)=1$. Now
the map
\begin{eqnarray*}
(\tilde{K}\times P)\times(\tilde{K}\times P)&\ra& \tilde{K}\times P\\
(k',p'),(k,p) &\mapsto&
(kk'\tilde{a}_K(k^{-1}p'k,p),a_P(k^{-1}p'k,p))
\end{eqnarray*}
defines a multiplication on $\tilde{G}=\tilde{K}\times P$ with
the desired properties.

Finally $\ker(\beta)$ will automatically be central because it is
a normal subgroup of order two.
\qed

Let $S$ be the spin representation of $Spin(\p)$ (see
\cite{LawMich}, p.36). It splits as a direct sum of two distinct
irreducible representations
$$
S\= S^+\oplus S^-.
$$
We will not go into the theory of the spin representation, we
only need the following properties:
\begin{itemize}
\item
The virtual representation
$$
(S^+- S^-)\otimes (S^+- S^-)
$$
is isomorphic to the adjoint representation on $\wedge^{even}\p
-\wedge^{odd}\p$ (see \cite{LawMich}, p. 36).
\item
If $K$ leaves invariant the spaces $\p_-$ and $\p_+$, as is the
case when $X$ carries a holomorphic structure fixed by $G$, then
there is a one dimensional representation $\epsilon$ of
$\tilde{K}$ such that
$$
(S^+-S^-)\otimes\epsilon \ \cong\
\wedge^{even}\p_--\wedge^{odd}\p_-.
$$
\end{itemize}
The proof of this latter property will be given in section
\ref{appA}.

\begin{theorem} \label{existh}
(Pseudo-coefficients) Assume that $G$ is orientation
preserving. Then for any finite dimensional representation
$(\tau ,V_\tau)$ of $\tilde{K}$ there is a compactly
supported smooth function $h_\tau$ on $\tilde{G}$ such
that for every irreducible unitary representation $(\pi
,V_\pi)$ of $\tilde{G}$ it holds:
$$
\tr\ \pi (h_\tau) = \dim (V_\pi \otimes S^+ \otimes
V_{\breve{\tau}})^{\tilde{K}} - \dim (V_\pi \otimes S^- \otimes
V_{\breve{\tau}})^{\tilde{K}}.
$$
\end{theorem}

The functions given in this theorem are also known as
pseudo-coefficients \cite{Lab}. This result generalizes
the one in \cite{Lab} in several ways. First, the group
$G$ needn't be connected and secondly the representation
$\tau$ needn't be spinorial. The proof of this theorem
relies on the following lemma.

\begin{lemma}
Let $(\pi ,V_\pi)$ be an irreducible unitary representation of
$\tilde{G}$ and assume
$$
\dim (V_\pi \otimes S^+ \otimes V_{\breve{\tau}})^{\tilde{K}} -
\dim (V_\pi \otimes S^- \otimes V_{\breve{\tau}})^{\tilde{K}}
\neq 0,
$$
then the Casimir eigenvalue satisfies $\pi (C) = \breve{\tau}
(C_K) - B(\rho)+B(\rho_K)$.
\end{lemma}

\prf Let the $\tilde{K}$-invariant operator
 $$
 d_\pm : V_\pi\otimes S^\pm \ra V_\pi\otimes S^\mp
 $$
be defined by
 $$
 d_\pm : v\otimes s \mapsto \sum_i\pi(X_i)v\otimes c(X_i)s,
 $$
where $(X_i)$ is an orthonormal base of $\p$. The formula of
Parthasarathy, \cite{AtSch}, p. 55 now says
 $$
 d_-d_+\= d_+d_- \= \pi\otimes s^\pm(C_K) -\pi(C)\otimes 1
 -B(\rho)+B(\rho_K).
 $$
Our assumption leads to $ker(d_+ d_-) \cap \pi \otimes S(\tau)
\neq 0$, and therefore $0=\tau(C_K)-\pi(C) -B(\rho) +B(\rho_K)$.
\qed

For the proof of Theorem \ref{existh} let $(\tau,V_\tau)$ a finite
dimensional irreducible unitary representation of $\tilde{K}$ and
write $E_\tau$ for the $\tilde{G}$-homogeneous vector bundle over
${X}=\tilde{G}/\tilde{K}$ defined by $\tau$. The space of smooth
sections $\Ga^\infty (E_\tau)$ may be written as $\Ga^\infty
(E_\tau) = (C^\infty (\tilde{G}) \otimes V_\tau)^{\tilde{K}}$,
where $\tilde{K}$ acts on $C^\infty (\tilde{G})$ by right
translations. The Casimir operator $C$ of $\tilde{G}$ acts on
this space and defines a second order differential operator
$C_\tau$ on $E_\tau$. On the space of $L^2$-sections
$L^2(X,E_\tau)=(L^2(\tilde{G}) \otimes V_\tau)^{\tilde{K}}$ this
operator is formally selfadjoint with domain, say, the compactly
supported smooth functions and extends to a selfadjoint operator.
Consider a Schwartz function $f$ on $\R$ such that the Fourier
transform $\hat{f}$ has compact support. By general results on
hyperbolic equations (\cite{Tayl}, chap IV) it follows that the
smoothing operator $f(C_\tau)$, defined by the spectral theorem
has finite propagation speed. Since $f(C_\tau)$ is
$\tilde{G}$-equivariant it follows that $f(C_\tau)$ can be
represented as a convolution operator $\ph \mapsto \ph *
\breve{f}_\tau$, for some $\breve{f}_\tau \in
(C_c^\infty(\tilde{G})\otimes \End (V_\tau))^{\tilde{K}\times
\tilde{K}}$ and with $\tilde{f}_\tau$ denoting the pointwise
trace of $\breve{f}_\tau$ we have for $\pi \in \hat{\tilde{G}}$:
$\tr\ \pi(\tilde{f}_\tau) = f(\pi(C))\dim(V_\pi \otimes
V_\tau)^{\tilde{K}}$, where $\pi(C)$ denotes the Casimir
eigenvalue on $\pi$. This construction extends to virtual
representations by linearity.

Choose $f$ such that $f(\breve{\tau} (C_K) -B(\rho)+B(\rho_K))=1$.
Such an $f$ clearly exists. Let $\ga$ be the virtual
representation of $\tilde{K}$ on the space
$$
V_\ga = (S^+-S^-)\otimes V_\tau,
$$
then set $h_\tau =\tilde{f}_\ga$. Theorem \ref{existh} follows.
\qed

To get the first part of Theorem \ref{exist-ep} from Theorem
\ref{existh} one replaces $\tau$ in the proposition by the virtual
representation on $(S^+-S^-)\otimes V_\tau$. Since
$(S^+-S^-)\otimes (S^+-S^-)$ is as $\tilde{K}$ module isomorphic
to $\wedge^*\p$ we get the desired function, say $j$ on the group
$\tilde{G}$. Now if $\tilde{G}\ne G$ let $z$ be the nontrivial
element in the kernel of the isogeny $\tilde{G}\ra G$, then the
function
$$
f(x) \= \rez{2}(j(x)+j(zx))
$$
factors over $G$ and satisfies the claim.

To get the second part of the theorem one proceeds similarly
replacing $\tau$ by $\epsilon\otimes \tau$.
\qed

\section{Clifford algebras and Spin groups}\label{appA}
This section is solely given to provide a proof of the properties
of the spin representation used in the last section. We will
therefore not strive for the utmost generality but plainly state
things in the form needed. For more details the reader is
referred to \cite{LawMich}.

Let $V$ be a finite dimensional complex vector space and let
$q:V\ra\C$ be a non-degenerate quadratic form. We use the same
letter for the symmetric bilinear form:
$$
q(x,y) \= \rez{2}(q(x+y)-q(x)-q(y)).
$$
Let $SO(q)\subset GL(V)$ be the special orthogonal group of $q$.
The {\it Clifford algebra} \index{Clifford algebra} $Cl(q)$ will
be the quotient of the tensorial algebra
$$
TV \= \C\oplus V\oplus (V\otimes V)\oplus\dots
$$
by the two-sided ideal generated by all elements of the form
$v\otimes v+q(v)$, where $v\in V$.

This ideal is not homogeneous with respect to the natural
$\Z$-grading of $TV$, but it is homogeneous with respect to the
induced $\Z/2\Z$-grading given by the even and odd degrees. Hence
the latter is inherited by $Cl(q)$:
$$
Cl(q)\= Cl^0(q)\oplus Cl^1(q).
$$
For any $v\in V$ we have in $Cl(q)$ that $v^2=-q(v)$ and therefore
$v$ is invertible in $Cl(q)$ if $q(v)\ne 0$. Let $Cl(q)^\times$
be the group of invertible elements in $Cl(q)$. The algebra
$Cl(q)$ has the following universal property: For any linear map
$\ph :V\ra A$ to a $\C$-algebra $A$ such that $\ph(v)^2=-q(v)$
for all $v\in V$ there is a unique algebra homomorphism $Cl(v)\ra
A$ extending $\ph$.

Let $Pin(q)$ be the subgroup of the group $Cl(q)^\times$
generated by all elements $v$ of $V$ with $q(v)=\pm 1$. Let the
{\it complex spin group} \index{complex spin group} be defined by
$$
Spin(q) \= Pin(q)\cap Cl^0(q),
$$
i.e., the subgroup of $Pin(q)$ of those elements which are
representable by an even number of factors of the form $v$ or
$v^{-1}$ with $v\in V$. Then $Spin(q)$ acts on $V$ by $x.v =
xvx^{-1}$ and this gives a double covering: $Spin(q)\ra SO(q)$.

Assume the dimension of $V$ is even and let
$$
V \= V^+\oplus V^-
$$
be a {\it polarization}, \index{polarization} that is
$q(V^+)=q(V^-)=0$. Over $\C$ polarizations always exist for even
dimensional spaces. By the nondegeneracy of $q$ it follows that
to any $v\in V^+$ there is a unique $\hat{v}\in V^-$ such that
$q(v,\hat{v})=-1$. Further, let $V^{-,v}$ be the space of all
$w\in V^-$ such that $q(v,w)=0$, then
$$
V^- \= \C \hat{v}\oplus V^{-,v}.
$$
Let
$$
S\= \wedge^* V^- \= \C\oplus V^-\oplus \wedge^2
V^-\oplus\dots\oplus\wedge^{top}V^-,
$$
then we define an action of $Cl(q)$ on $S$ in the following way:
\begin{itemize}
\item
for $v\in V^-$ and $s\in S$ let
$$
v.s\= v\wedge s,
$$
\item
for $v\in V^+$ and $s\in \wedge^*V^{-,v}$ let
$$
v.s\= 0,
$$
\item
and for $v\in V^+$ and $s\in S$ of the form $s=\hat{v}\wedge s'$
with $s'\in \wedge^*V^{-,v}$ let
$$
v.s\= s'.
$$
\end{itemize}
By the universal property of $Cl(V)$ this extends to an action of
$Cl(q)$. The module $S$ is called the {\it spin
module}.\index{spin module} The induced action of $Spin(q)$
leaves invariant the subspaces
$$
S^+\=\wedge^{even}V^-,\ \ \ S^-\=\wedge^{odd}V^-,
$$
the representation of $Spin(q)$ on these spaces are called the
{\it half spin representations}.\index{half spin representations}
Let $SO(q)^+$ the subgroup of all elements in $SO(q)$ that leave
stable the decomposition $V=V^+\oplus V^-$. This is a connected
reductive group isomorphic to $GL(V^+)$, since, let $g\in
GL(V^+)$ and define $\hat{g}\in GL(V^-)$ to be the inverse of the
transpose of $g$ by the pairing induced by $q$ then the map
$Gl(v)\ra SO(q)^+$ given by $g\mapsto (g,\hat{g})$ is an
isomorphism. In other words, choosing a basis on $V^+$ and a the
dual basis on $V^-$ we get that $q$ is given in that basis by
$\matrix{0}{\1}{\1}{0}$. Then $SO(q)^+$ is the image of the
embedding
\begin{eqnarray*}
GL(V^-) &\hookrightarrow& SO(q)\\
A &\mapsto& \matrix{A}{0}{0}{^tA^{-1}}.
\end{eqnarray*}

Let $Spin(q)^+$ be the inverse image of $SO(q)^+$ in $Spin(q)$.
Then the covering $Spin(q)^+\ra SO(q)^+\cong GL(V^-)$ is the
``square root of the determinant'', i.e., it is isomorphic to the
covering $\tilde{GL}(V^-)\ra GL(V^-)$ given by the pullback
diagram of linear algebraic groups:
$$
\begin{array}{ccc}
\tilde{GL}(V^-)     & \ra   & GL(1)\\
\downarrow         & {} &\downarrow x\mapsto x^2\\
GL(V^-)         & \begin{array}{c}\det\\ \ra\\ {}\end{array}&
GL(1).
\end{array}
$$
As a set, $\tilde{GL}(V^-)$ is given  as the set of all pairs
$(g,z)\in GL(V^-)\times GL(1)$ such that $\det(g)=z^2$ and the
maps to $GL(V^-)$ and $GL(1)$ are the respective projections.

\begin{lemma}\label{epsilon}
There is a one dimensional representation $\epsilon$ of
$Spin(q)^+$ such that
$$
S^\pm\otimes\epsilon\ \cong\ \wedge^\pm V
$$
as $Spin(q)^+$-modules, where $\wedge^\pm$ means the even or odd
powers respectively.
\end{lemma}

\prf Since $Spin(q)^+$ is a connected reductive group over $\C$
we can apply highest weight theory. If the weights of the
representation of $Spin(q)^+$ on $V$ are given by
$\pm\mu_1,\dots,\pm\mu_m$, then the weights of the half spin
representations are given by
$$
\rez{2}(\pm\mu_1\pm\dots\pm\mu_m)
$$
with an even number of minus signs in the one and an odd number
in the other case. Let $\epsilon = \rez{2}(\mu_1+\dots +\mu_m)$
then $\epsilon$ is a weight for $Spin(q)^+$ and $2\epsilon$ is
the weight of, say, the one dimensional representation on
$\wedge^{top}V^+$. By Weyl's dimension formula this means that
$2\epsilon$ is invariant under the Weyl group and therefore
$\epsilon$ is. Again by Weyl's dimension formula it follows that
the representation with highest weight $\epsilon$ is one
dimensional. Now it follows that $S^+\otimes\epsilon$ has the
same weights as the representation on $\wedge^+ V$, hence must be
isomorphic to the latter. The case of the minus sign is analogous.
\qed

\section{Orbital integrals}
It now will be shown that $\tr\pi(f_\tau)$ vanishes for a
principal series representation $\pi$. To this end let $P=MAN$ be
a nontrivial parabolic subgroup with $A\subset \exp(\p_0)$. Let
$(\xi ,V_\xi)$ denote an irreducible unitary representation of $M$
and $e^\nu$ a quasicharacter of $A$. Let $\pi_{\xi ,\nu}:= {\rm
Ind}_P^G (\xi \otimes e^{\nu}\otimes 1)$.

\begin{lemma} \label{pivonggleichnull}
We have $\tr\pi_{\xi ,\nu}(f_\tau) =0$.
\end{lemma}

\prf By Frobenius reciprocity we have for any irreducible unitary
representation $\ga$ of $K$:
$$
\Hom_K(\ga ,\pi_{\xi ,\nu}|_K) \cong \Hom_{K_M}(\ga |_{K_M},\xi ),
$$
where $K_M := K\cap M$. This implies that $\tr\pi_{\xi
,\nu}(f_\tau)$ does not depend on $\nu$. On the other hand
$\tr\pi_{\xi ,\nu}(f_\tau)\ne 0$ for some $\nu$ would imply
$\pi_{\xi ,\nu}(C) =\breve{\tau} (C_K) -B(\rho)+B(\rho_K)$ which
only can hold for $\nu$ in a set of measure zero.
\qed

Recall that an element $g$ of $G$ is called {\it
elliptic}\index{elliptic} if it lies in a compact Cartan subgroup.
Since the following relies on results of Harish-Chandra which were
proven under the assumption that $G$ is of inner type, we will
from now on assume this.

\begin{proposition} \label{orbitalint}
Assume that $G$ is of inner type. Let $g$ be a semisimple element
of the group $G$. If $g$ is not elliptic, then the orbital
integral $\O_g(f_\tau)$ vanishes. If $g$ is elliptic we may assume
$g\in T$, where $T$ is a Cartan in $K$ and then we have
$$
\O_g(f_\tau) \= \tr\ \tau(g)\ c_g^{-1}|W(\t ,\g_g)| \prod_{\alpha
\in \Phi_g^+}(\rho_g ,\alpha),
$$
where $c_g$ is Harish-Chandra's constant, it does only depend on
the centralizer $G_g$ of g. Its value is given in \ref{notations}.

\end{proposition}

\prf The vanishing of $\O_g(f_\tau)$ for nonelliptic semisimple
$g$ is immediate by the lemma above and Theorem \ref{inv-orb-int}.
So consider $g\in T\cap G'$, where $G'$ denotes the set of regular
elements. Note that for regular $g$ the claim is $\CO_g(f_\tau)
=\tr \tau(g)$. Assume the claim proven for regular elements, then
the general result follows by standard considerations as in
\cite{HC-DS}, p.32 ff. where however different Haar-measure
normalizations are used that produce a factor $[G_g:G_g^0]$,
therefore these standard considerations are now explained. Fix
$g\in T$ not necessarily regular. Let $y\in T^0$ be such that $gy$
is regular. Then
\begin{eqnarray*}
\tr \tau (gy) &=& \int_{T\bs G} f_\tau(x^{-1}gyx) dx\\
    &=&  \int_{T^0\bs G} f_\tau(x^{-1}gyx) dx\\
    &=&  \int_{G_g\bs G}\int_{T^0\bs G_g} f_\tau(x^{-1}z^{-1}gyzx)\ dz\ dx\\
    &=&  \int_{G_g\bs G}\sum_{\eta :G_g /G_g^0}
\rez{[G_g:G_g^0]}\int_{T^0\bs G_g^0}
f_\tau(x^{-1}\eta^{-1}z^{-1}gyz\eta x)\ dz\ dx.
\end{eqnarray*}
The factor $\rez{[G_g:G_g^0]}$ comes in by the Haar-measure
normalizations. On $G_g^0$ consider the function
$$
h(y) = f(x^{-1}\eta^{-1} yg\eta x).
$$
Now apply Harish-Chandra's operator $\omega_{G_g}$ to $h$ then
for the connected group $G_g^0$ it holds
$$
h(1) = \lim_{y\ra 1} c_{G_g^0}\ \omega_{G_g^0}\ \CO_y^{G_g^0}(h).
$$
When $y$ tends to $1$ the $\eta$-conjugation drops out and the
claim follows.

So in order to prove the proposition one only has to consider the
regular orbital integrals. Next the proof will be reduced to the
case when the compact Cartan $T$ meets all connected components
of $G$. For this let $G^+=TG^0$ and assume the claim proven for
$G^+$. Let $x\in G$ then $xTx^{-1}$ again is a compact Cartan
subgroup. Since $G^0$ acts transitively on all compact Cartan
subalgebras it follows that $G^0$ acts transitively on the set of
all compact Cartan subgroups of $G$. It follows that there is a
$y\in G^0$ such that $xTx^{-1} =yTy^{-1}\subset TG^0 =G^+$, which
implies that $G^+$ is normal in $G$.

Let $\tau^+=\tau |_{G^+\cap K}$ and $f^+_{\tau^+}$ the
corresponding Euler-Poincar\'e function on $G^+$.

\begin{lemma}
$f^+_{\tau^+} = f_\tau |_{G^+}$
\end{lemma}

Since the Euler-Poincar\'e function is not uniquely determined the
claim reads that the right hand side is a EP-function for $G^+$.

\prf Let $\tau^+ =\tau |_{K^+}$, where $K^+ = TK^0 = K\cap G^+$.
Let $\ph^+\in(C_c^\infty(G^+)\otimes V_\tau)^{K^+}$, which may be
viewed as a function $\ph^+ : G^+ \ra V_\tau$ with $\ph^+(xk)
=\tau(k^{-1})\ph^+(x)$ for $x\in G^+$, $k\in K^+$. Extend $\ph^+$
to $\ph : G\ra V_\tau$ by $\ph(xk)=\tau(k^{-1})\ph^+(x)$ for
$x\in G^+$, $k\in K$. This defines an element of
$(C_c^\infty(G)\otimes V_\tau)^K$ with $\ph |_{G^+}=\ph^+$. Since
$C_\tau$ is a differential operator it follows $f(C_\tau)\ph
|_{G^+} =f(C_\tau)\ph^+$, so
$$
(\ph * \tilde{f}_\tau) |_{G^+} = \ph^+ * \tilde{f}_{\tau^+}.
$$
Considering the normalizations of Haar measures gives the lemma.
\qed

For $g\in T'$ we compute
\begin{eqnarray*}
\CO_g(f_\tau) &=& \int_{T\bs G} f_\tau (x^{-1} gx) dx\\
    &=& \sum_{y:G/G^+} \rez{[G:G^+]} \int_{T\bs G^+} f_\tau (x^{-1} y^{-1} g yx) dx,
\end{eqnarray*}
 where the factor $\rez{[G:G^+]}$ stems from normalization of
Haar measures and we have used the fact that $G^+$ is normal. The
latter equals
$$
\rez{[G:G^+]}\sum_{y:G/G^+}\CO^{G^+}_{y^{-1}xy}(f_\tau) =
\rez{[G:G^+]}\sum_{y:G/G^+}\CO^{G^+}_{y^{-1}xy}(f_{\tau^+}^+).
$$
Assuming the proposition proven for $G^+$, this is
$$
\rez{[G:G^+]}\sum_{y:G/G^+} \tr \tau(y^{-1}gy) =\tr\tau(g).
$$

From now on one thus may assume that the compact Cartan $T$ meets
all connected components of $G$. Let $(\pi,V_\pi)\in\hat{G}$.
Harish-Chandra has shown that for any $\ph\in C^\infty_c(G)$ the
operator $\pi(\ph)$ is of trace class and there is a locally
integrable conjugation invariant function $\Theta_\pi$ on $G$,
smooth on the regular set such that
$$
\tr \pi(\ph) = \int_G\ph(x)\Theta_\pi(x) dx.
$$
For any $\psi\in C^\infty(K)$ let
$\pi|_K(\psi)=\int_K\psi(k)\pi(k) dk$.

\begin{lemma} \label{character}
Assume $T$ meets all components of $G$. For any $\psi\in
C^\infty(K)$ the operator $\pi|_K(\psi)$ is of trace class and
for $\psi$ supported in the regular set $K' =K\cap G'$ we have
$$
\tr \pi|_K(\psi) = \int_K \psi(k) \Theta_\pi(k) dk.
$$
(For $G$ connected this assertion is in \cite{AtSch} p.16.)
\end{lemma}

\prf Let $V_\pi =\bigoplus_i V_\pi(i)$ be the decomposition of
$V_\pi$ into $K$-types. This is stable under $\pi|_K(\psi)$.
Harish-Chandra has proven $[\pi|_K :\tau]\le\dim \tau$ for any
$\tau\in\hat{K}$. Let $\psi =\sum_j\psi_j$ be the decomposition
of $\psi$ into $K$-bitypes. Since $\psi$ is smooth the sequence
$\norm{\psi_j}_1$ is rapidly decreasing for any enumeration of
the $K$-bitypes. Here $\norm{\psi}_1$ is the $L^1$-norm on $K$.
It follows that the sum $\sum_i\tr(\pi|_K(\psi)|V_\pi(i))$
converges absolutely, hence $\pi|_K(\psi)$ is of trace class.

Now let $S=\exp(\p_0)$ then $S$ is a smooth set of
representatives of $G/K$. Let $G.K=\cup_{g\in
G}gKg^{-1}=\cup_{s\in S}sKs^{-1}$, then, since $G$ has a compact
Cartan, the set $G.K$ has non-empty interior. Applying the Weyl
integration formula to $G$ and backwards to $K$ gives the
existence of a smooth measure $\mu$ on $S$ and a function $D$ with
$D(k)>0$ on the regular set such that
$$
\int_{G.K}\ph(x) dx = \int_S \int_K \ph(sks^{-1}) D(k) dk d\mu(s)
$$
for $\ph\in L^1(G.K)$. Now suppose $\ph\in C_c^\infty(G)$ with
support in the regular set. Then
\begin{eqnarray*}
\tr \pi(\ph) &=& \int_{G.K}\ph(x) \Theta_\pi(x) dx\\
    &=& \int_S \int_K \ph(sks^{-1}) D(k) \Theta_\pi(k)d\mu(s)\\
    &=& \int_K\int_S \ph^s(k) d\mu(s) D(k) \Theta_\pi(k) dk,
\end{eqnarray*}
where we have written $\ph^s(k)=\ph(sks^{-1})$. On the other hand
\begin{eqnarray*}
\tr\pi(\ph) &=& \tr \int_{G.K} \ph(x) \pi(x) dx\\
    &=& \tr \int_S\int_K \ph(sks^{-1})D(k) \pi(sks^{-1})dk d\mu(s)\\
    &=& \tr \int_S \pi(s) \pi|_K(\ph^s D)\pi(s)^{-1} d\mu(s)\\
    &=& \int_S \tr \pi |_K(\ph^s D) d\mu(s)\\
    &=& \tr\pi |_K\left( \int_S \ph^s d\mu(s) D\right).
\end{eqnarray*}
This implies the claim for all functions $\psi\in C_c^\infty(K)$
which are of the form
$$
\psi(k) = \int_S\ph(sks^{-1})d\mu(s) D(k)
$$
for some $\ph\in C_c^\infty(G)$ with support in  the regular set.
Consider the map
$$
\begin{array}{cccc}
F:& S\times K' & \ra & G.K'\\
{}& (s,k) & \mapsto & sks^{-1}
\end{array}
$$
Then the differential of $F$ is an isomorphism at any point and by
the inverse function theorem $F$ locally is a diffeomorphism. So
let $U\subset S$ and $W\subset K'$ be open sets such that
$F|_{U\times W}$ is a diffeomorphism. Then let $\alpha\in
C_c^\infty(U)$ and $\beta\in C_c^\infty(W)$, then define
$$
\phi(sks^{-1}) = \alpha(s)\beta(k)\ \ \ {\rm if}\ s\in U,\ k\in W
$$
and $\ph(g)=0$ if $g$ is not in $F(U\times W)$. We can choose the
function $\alpha$ such that $\int_S\alpha (s) d\mu(s) =1$. Then
$$
\int_S\ph(sks^{-1}) d\mu(s) D(k) = \beta(k) D(k).
$$
Since $\beta$ was arbitrary and $D(k) >0$ on $K'$ the lemma
follows.
\qed

Let $W$ denote the virtual $K$-representation on $\wedge^{even}\p
\otimes V_\tau - \wedge^{odd}\p\otimes V_\tau$ and write $\chi_W$
for its character.

\begin{lemma} \label{fin-lin-comb}
Assume $T$ meets all components of $G$, then for any
$\pi\in\hat{G}$ the function $\Theta_\pi\chi_W$ on $K'=K\cap G'$
equals a finite integer linear combination of $K$-characters.
\end{lemma}

\prf It suffices to show the assertion for $\tau =1$. Let $\ph$
be the homomorphism $K\ra O(\p)$ induced by the adjoint
representation, where the orthogonal group is formed with respect
to the Killing form. We claim that $\ph(K)\subset SO(\p)$, the
subgroup of elements of determinant one. Since we assume $K=K^0
T$ it suffices to show $\ph(T)\subset SO(\p)$. For this let $t\in
T$. Since $t$ centralizes $\t$ it fixes the  decomposition $\p
=\oplus_\alpha \p_\alpha$ into one dimensional root spaces. So
$t$ acts by a scalar, say $c$ on $\p_\alpha$ and by $d$ on
$\p_{-\alpha}$. There is $X\in\p_\alpha$ and $Y\in\p_{-\alpha}$
such that $B(X,Y)=1$. By the invariance of the Killing form $B$
we get
$$
1 = B(X,Y) = B(\Ad(t)X,\Ad(t)Y) = cdB(X,Y) =cd.
$$
So on each pair of root spaces $\Ad(t)$ has determinant one hence
also on $\p$.

Replacing $G$ by a double cover if necessary, which doesn't effect
the claim of the lemma, we may assume that $\ph$ lifts to the spin
group $\Spin(\p)$. Let $\p=\p^+\oplus\p^-$ be the decomposition
according to an ordering of $\phi(\t ,\g)$. This decomposition is
a polarization of the quadratic space $\p$ and hence the spin
group acts on $S^+=\wedge^{even}\p^+$ and $S^-=\wedge^{odd}\p^+$
in a way that the virtual module $(S^+-S^-)\otimes (S^+-S^-)$
becomes isomorphic to $W$. For $K$ connected the claim now follows
from \cite{AtSch} (4.5). An inspection shows however that the
proof of (4.5) in \cite{AtSch}, which is located in the appendix
(A.12), already applies when we only assume that the homomorphism
$\ph$ factors over the spin group. \qed

We continue the proof of the proposition. Let $\hat{T}$ denote the
set of all unitary characters of $T$. Any regular element
$\hat{t}\in \hat{T}$ gives rise to a discrete series
representation $(\omega ,V_\omega)$ of $G$. Let
$\Theta_{\hat{t}}=\Theta_\omega$ be its character which, due to
Harish-Chandra, is known to be a function on $G$. Harish-Chandra's
construction gives a bijection between the set of discrete series
representations of $G$ and the set of $W(G,T)=W(K,T)$-orbits of
regular characters of $T$.

Let $\Phi^+$ denote the set of positive roots of $(\g,\t)$ and let
$\Phi^+_c$ ,$\Phi^+_n$ denote the subsets of compact and
noncompact positive roots. For each root $\alpha$ let $t\mapsto
t^\alpha$ denote the corresponding character
 on $T$.
Define
\begin{eqnarray*}
'\Delta_c(t) &=& \prod_{\alpha\in\Phi^+_c}(1-t^{-\alpha})\\
'\Delta_n &=& \prod_{\alpha\in\Phi^+_n}(1-t^{-\alpha})
\end{eqnarray*}
and $'\Delta = '\Delta_c '\Delta_n$. If $\hat{t}\in \hat{T}$ is
singular, Harish-Chandra has also constructed an invariant
distribution $\Theta_{\hat{t}}$ which is a virtual character on
$G$. For $\hat{t}$ singular let $W(\hat{t})\subset W(\g,\t)$ be
the isotropy group. One has $\Theta_{\hat{t}} = \sum_{w\in
W(\hat{t})} \epsilon(w)\Theta'_{w,\hat{t}}$ with
$\Theta'_{w,\hat{t}}$ the character of an induced representation
acting on some Hilbert space $V_{w,\hat{t}}$ and $\epsilon(w)\in
\{\pm 1\}$. Let $\CE_2(G)$ denote the set of discrete series
representations of $G$ and $\CE_2^s(G)$ the set of $W(G,T)$-orbits
of singular characters.

By Theorem \ref{inv-orb-int} the proposition will follow from the

\begin{lemma}\label{trace-tau}
For $t\in T$ regular we have
$$
\tr \tau(t)\=
\rez{|W(G,T)|}\sum_{\hat{t}\in\hat{T}}\Theta_{\hat{t}}(f_\tau)
\Theta_{\hat{t}}(t).
$$
\end{lemma}

\prf Let $\ga$ denote the virtual $K$-representation on
$(\wedge^{even}\p -\wedge^{odd}\p)\otimes V_\tau$. Harish-Chandra
has shown (\cite{HC-S} Theorem 12) that for any
$\hat{t}\in\hat{T}$ there is an irreducible unitary representation
$\pi_{\hat{t}}^0$ such that $\Theta_{\hat{t}}$ coincides up to
sign with the character of $\pi_{\hat{t}}^0$ on the set of
elliptic elements of $G$ and $\pi_{\hat{t}}^0=\pi_{\hat{t}'}^0$ if
and only if there is a $w\in W(G,T)=W(K,T)$ such that
$\hat{t}'=w\hat{t}$.

Further (\cite{HC-S}, Theorem 14) Harish-Chandra has shown that
the family
$$
\left( \frac{'\lap(t)\Theta_{\hat{t}}(t)}
            {\sqrt{|W(G,T)|}}\right)_{\hat{t}\in \hat{T}/W(G,T)}
$$
forms an orthonormal basis of $L^2(T)$. Here we identify
$\hat{T}/W(G,T)$ to a set of representatives in $\hat{T}$ to make
$\Theta_{\hat{t}}$ well defined.

Consider the function $g(t)=\frac{\tr\ga(t) \
'\lap_c(t)}{\overline{\ '\lap_n(t)}} = \tr\tau(t) \ '\lap(t)$. Its
coefficients with respect to the above orthonormal basis are
\begin{eqnarray*}
\langle g, \frac{\ '\lap\Theta_{\hat{t}}}
            {\sqrt{|W(G,T)|}}\rangle &=&
            \rez{\sqrt{|W(G,T)|}}\int_T\tr\ga(t) |
            '\lap_c(t)|^2\overline{\Theta_{\hat{t}}(t)} dt\\
 &=& \sqrt{|W(G,T)|}\int_K \tr\ga(k)\overline{\Theta_{\hat{t}}(k)}
 dk
\end{eqnarray*}
where we have used the Weyl integration formula for the group $K$
and the fact that $W(G,T)=W(K,T)$. Next by Lemma
\ref{fin-lin-comb} this equals
$$
\sqrt{|W(G,T)|}
\dim((\wedge^{even}\p-\wedge^{odd}\p)\otimes\breve{\tau}\otimes\pi_{\hat{t}}^0)^K
\= \sqrt{|W(G,T)|}\Theta_{\hat{t}}(f_\tau).
$$
Hence
\begin{eqnarray*}
g(t) &=& \tr\tau(t) \ '\lap(t)\\
    &=& \sum_{\hat{t}\in \hat{T}/W(G,T)}\Theta_{\hat{t}}(f_\tau)
    \ '\lap(t) \Theta_{\hat{t}}(t)\\
    &=& \rez{|W(G,T)|}\sum_{\hat{t}\in \hat{T}}\Theta_{\hat{t}}(f_\tau)
    \ '\lap(t) \Theta_{\hat{t}}(t).
\end{eqnarray*}
The lemma and the proposition are proven. \qed

\begin{corollary}
If $\tilde{g}\in\tilde{G}$ is semisimple and not elliptic then
$\CO_{\tilde{g}}(g_\tau)=0$. If $\tilde{g}$ is elliptic regular
then
$$
\CO_{\tilde{g}}(g_\tau)\=\frac{\tr\tau(\tilde{g})} {\tr(\tilde{g}
| S^+-S^-)}.
$$
\end{corollary}

\prf Same as for the last proposition with $g_\tau$ replacing
$f_\tau$.
\qed

\begin{proposition}
Assume that $\tau$ extends to a representation of the group $G$
on the same space. For the function $f_{\tau}$ we have for any
$\pi \in \hat{G}$:
$$
\tr \ \pi(f_{{\tau}}) \= \sum_{p=0}^{\dim \ \g /\k}(-1)^p \dim \
{\rm Ext}_{(\g ,K)}^p (V_{{\tau}} ,V_\pi),
$$
i.e., $f_{{\tau}}$ gives the Euler-Poincar\'e numbers of the $(\g
,K)$-modules $(V_{{\tau}} ,V_\pi)$, this justifies the name
Euler-Poincar\'e function.
\end{proposition}

\prf By definition it is clear that
$$
\tr \ \pi (f_{{\tau}}) \= \sum_{p=0}^{\dim \ \p} (-1)^p \dim \
H^p(\g ,K,V_{\breve{\tau}} \otimes V_\pi).
$$
The claim now follows from \cite{BorWall}, p. 52. \qed

%\hglhg
\end{document}